\documentclass[preprint,11pt]{elsarticle}

\newtheorem{theorem}{Theorem}[section]

\newtheorem{lemma}{Lemma}[section]

\usepackage{amssymb}
\usepackage{amsmath}

\begin{document}

\begin{frontmatter}

\title{Bernstein type inequalities for   self-normalized martingales  with applications}
\author{Xiequan Fan$^*$}
\author{  Shen Wang}
 \cortext[cor1]{\noindent Corresponding author. \\
  \mbox{ \ \ \ } \textit{E-mail}: fanxiequan@hotmail.com (X. Fan). }
\address{Center for Applied Mathematics, Tianjin University, 300072 Tianjin,  China}

\begin{abstract}
For self-normalized martingales with conditionally symmetric differences, de la Pe\~{n}a \cite{VH99} established the Gaussian type exponential inequalities. Bercu and Touati \cite{BT08}   extended  de la Pe\~{n}a's inequalities to martingales with differences heavy on left.
In this paper, we establish Bernstein type exponential inequalities for self-normalized martingales with differences bounded from below.
 Moreover, applications to  self-normalized sums, t-statistics and autoregressive processes  are discussed.
\end{abstract}

\begin{keyword}  martingales;  self-normalized processes;  exponential inequalities; autoregressive processes
\vspace{0.3cm}
\MSC Primary 60G42; 60E15; 60F10
\end{keyword}

\end{frontmatter}


\section{Introduction}
Let $ (\xi_i)_{i \geq1}$ be a sequence of zero-mean independent random
variables satisfying $ \xi_i \leq 1 $   for
all $i$. Denote $S_n=\sum_{i=1}^{n} \xi_i$   the partial sums of $ (\xi_i)_{i \geq1}$.  Bennett \cite{B62} proved the following Bernstein type inequality: for all $x>0,$
\begin{eqnarray} \label{sdfsdkm}
 \mathbb{P}\big(S_n \geq x v^2  \big) \leq     \exp \bigg\{-\frac{x^2 v^2}{ 2(1+ x  /3  )}\bigg\},
\end{eqnarray}
where $v^2=\textrm{Var}(S_n)$ is the variance of $S_n.$  The importance of Bernstein type inequalities comes from
the fact that they combine both the Gaussian trends and exponentially decaying rate. To see this,  we rewrite the last inequality  in the following form: for all $x>0,$
\begin{eqnarray}\label{sdfsdkn}
 \mathbb{P}\big(S_n \geq x    \big) \leq     \exp \bigg\{-\frac{x^2 }{ 2(v^2+ x  / 3    )}\bigg\} .
\end{eqnarray}
It is easy to see that the last bound behaves as $\exp\{- \frac{x^2}{2v^2} \}$  for moderate $x=o(v^2),$  while  it is exponentially decaying to $0$
as   $x\rightarrow \infty$.

The generalizations of (\ref{sdfsdkm}) to martingales have attracted certain interest.  Assume that $(\xi_i,\mathcal{F}_{i})_{i=0,\cdots,n}$ is a sequence of   martingale differences. If $\xi_i \leq 1,$ Freedman \cite{Fr75} showed that  (\ref{sdfsdkm}) holds also when $\mathbb{P}\big(S_n \geq x v^2 \big)$ is replaced by $\mathbb{P}\big(S_n \geq x v^2, \langle S\rangle_n \leq v^2  \big),$ where $\langle S\rangle_n$ is the conditional variance of $S_n.$ De la Pe\~{n}a \cite{VH99},   Dzhaparidze and van Zanten \cite{Dz01}  and Fan et al.\ \cite{FGL12,FGL17}   extended Freedman's inequality to martingales with non-bounded differences.  Recently, Rio \cite{R17} gave a refinements on Freedman's inequality.

Despite the fact that the case of martingale is well studied, there are only a few results on Bernstein type inequalities for  self-normalized  martingales $ S_n/ [S]_n $, where $[S]_n$ is the squared variance of $S_n$.  Among them, let us recall the following exponential inequalities of de la Pe\~{n}a \cite{VH99}.
Assume that $(\xi_i,\mathcal{F}_{i})_{i=0,\cdots,n}$ is a sequence of conditionally symmetric martingale differences. Recall that $\xi_i$ is called
conditionally symmetric if $ \mathcal{L}(\xi_i|\mathcal{F}_{i-1}) = \mathcal{L}( -\xi_i|\mathcal{F}_{i-1})$ for all $i$,
where $\mathcal{L}(\xi_i|\mathcal{F}_{i-1})$ stands for  the regular version of the conditional distribution of $\xi_i$
given a $\sigma$-field $\mathcal{F}_{i-1}$.
De la Pe\~{n}a \cite{VH99}  have established   the following  exponential inequalities for self-normalized martingales:  for all $x> 0$,
\begin{equation}\label{ineqBT1}
\mathbb{P}\bigg( \frac{S_n}{  [S]_n }\geq x  \bigg) \leq  \sqrt{\mathbb{E} \bigg[   \exp{\bigg\{ -\frac12 x^2[S]_n \bigg\}}\bigg]} ,
\end{equation}
and, for all $x, y> 0$,
\begin{equation}\label{ineqBT2}
\mathbb{P}\bigg( \frac{S_n}{  [S]_n }\geq x, [S]_n \geq y  \bigg) \leq       \exp{\bigg\{ -\frac12x^2 y  \bigg\}}  ,
\end{equation}
where $[S]_n=\sum_{i=1}^n \xi _i^2$ the squared variance of $S_n.$
 In the i.i.d.\ case,  $[S]_n / n$ usually converges almost surely to the variance of the random variables. Thus (\ref{ineqBT1}) and (\ref{ineqBT2}) can be regarded  as  Gaussian type inequalities.

The  inequalities of de la Pe\~{n}a have been extended to the martingales with differences heavy on left.
Recall that an integrable random variable $X$  is called heavy on left if $\mathbb{E} X  = 0$ and, for all $a >0,$ $\mathbb{E}[T_a(X)] \leq 0,$
where
$$T_a(X) = \min(|X|, a) \,  \textrm{sign}(X) $$
is the truncated version of $X$.  Clearly,  conditionally symmetric martingale differences are heavy on left.
Bercu and Touati \cite{BT08} have obtained the following extension of de la Pe\~{n}a's inequity  (\ref{ineqBT1}):    for all $x> 0$,
\begin{equation}\label{ineqBT3}
\mathbb{P}\bigg( \frac{S_n}{  [S]_n }\geq x  \bigg) \leq \inf_{p >1} \bigg( \mathbb{E} \bigg[   \exp{\bigg\{ -\frac12(p-1)x^2[S]_n \bigg\}}\bigg] \bigg)^{  1 /p}.
\end{equation}
They also showed that (\ref{ineqBT2}) holds for martingales with differences heavy on left. In the particular case $p=2,$  inequality (\ref{ineqBT3}) reduces to
inequality (\ref{ineqBT1})  under the conditional symmetric assumption. Similar results for self-normalized martingales $S_n/\sqrt{[S]_n}$ can also
be found in  Bercu and Touati \cite{BT08}.


Exponential inequalities for  self-normalized martingales have a lot of applications.
We refer to \mbox{de la Pe\~{n}a,}  Klass and  Lai \cite{DML07} for autoregressive  processes. 
Bercu and Touati \cite{BT08} applied such type inequalities to parameter estimations of linear regressions, autoregressive  processes and branching processes. For more applications of such type inequalities, we refer to   the  monographs of   de la  Pe\~{n}a,  Lai and Shao \cite{DLS08}  and Bercu, Delyon and Rio \cite{BDR15}.

 In this paper, we aim to establish  Bernstein type inequalities for  self-normalized  martingales with differences bounded from below.
 It is obvious that a random variable is bounded from below does not imply that it is heavy on left.
Our results for self-normalized martingales are analogues  to the inequalities (\ref{ineqBT1}) - (\ref{ineqBT3}).
Applications to  self-normalized  sums,  t-statistics  and autoregressive processes are also discussed.

The paper is organized as follows. We present our main results in Section 2. In Section 3, we discuss the applications, and prove our main results in Section 4.

\section{Main results}
Let $(\xi_i,\mathcal{F}_{i})_{i=0,\cdots,n}$  be  a finite sequence of real-valued  square integrable martingale differences defined on a probability space $(\Omega,\mathcal{F},P)$, where $\xi_0=0$ and $\{\emptyset,\Omega\}=\mathcal{F}_{0}\subseteq\ldots\subseteq\mathcal{F}_{n}\subseteq\mathcal{F}$ are increasing $\sigma$-fields. So by definition, we have $\mathbb{E}[\xi_i|\mathcal{F}_{i-1}]=0,i=1,\ldots,n$. Set
\begin{equation}
S_0=0 \quad \textrm{and} \quad S_k=\sum_{i=1}^k\xi_i
\end{equation}
for $k=1,\ldots,n.$ Then $S=(S_k,\mathcal{F}_{k})_{k=1,\ldots,n}$ is a martingale.
Let $[ S]$ and $\langle S \rangle$  be, respectively, the squared variance and the conditional variance  of the
martingale $S,$ that is
\begin{eqnarray*}
[S]_0=0,\ \ \ \ \ [ S]_k=\sum_{i=1}^k \xi _i^2,
\end{eqnarray*}
and
\begin{equation}\label{quad01}
\langle S\rangle_0=0,\ \ \ \ \ \langle S \rangle_k=\sum_{i=1}^k \mathbb{E} [ \xi _i^2  |  \mathcal{F}_{i-1} ]  ,\quad k=1,...,n.
\end{equation}

Our main result is the following Bernstein type   inequalities for  self-normalized   martingales  with differences bounded from below. It is worth to  be mentioned that the inequalities
are new even for   independent random variables.
\begin{theorem} \label{thPSn}
Assume that $\xi_i \geq -1$ for all $i \in [1, n]$. Then for all $x>0$,
\begin{eqnarray}
\mathbb{P}\bigg( \frac{S_n}{[S]_n} \geq x \bigg)&\leq & \inf_{p>1}\bigg(\mathbb{E}\bigg[ \exp\bigg\{-(p-1) \Big( x-\log(1+x) \Big ) [S]_n  \bigg\} \mathbf{1}_{\{S_n   \geq x[S]_n\} } \bigg]  \bigg)^{ 1 /p}  \\
&\leq & \inf_{p>1}\bigg(\mathbb{E}\bigg[ \exp\bigg\{-(p-1)\frac{x^2}{2(1+x)} [S]_n \bigg\} \mathbf{1}_{\{S_n   \geq x[S]_n\} } \bigg]  \bigg)^{ 1 /p}, \label{ineqPSn01}
\end{eqnarray}
and, for all $  y> 0$,
\begin{eqnarray}
\mathbb{P} \bigg(\frac{S_n}{[S]_n} \geq x, \     [S]_n \geq y \bigg) &\leq & \exp\bigg\{- ( x-\log(1+x) ) y  \bigg\}   \label{ineqPSn012}\\
&\leq & \exp\bigg\{-\frac{x^2 y}{2(1+x)}  \bigg\}.\label{ineqPSn02}
\end{eqnarray}
\end{theorem}

Clearly, inequality (\ref{ineqPSn01}) implies that for all $x>0$,
 \begin{eqnarray*}
\mathbb{P}\bigg( \frac{S_n}{[S]_n} \geq x \bigg)&\leq & \inf_{p>1}\bigg(\mathbb{E}\bigg[ \exp\bigg\{-(p-1)\frac{x^2}{2(1+x)} [S]_n \bigg\} \bigg]  \bigg)^{ 1 /p},
\end{eqnarray*}
which is an analogues to de la Pe\~{n}a's inequity  (\ref{ineqBT1})  and the inequality of Bercu and Touati (\ref{ineqBT3}).

Denote $B_n^2= \sum_{i=1}^n\mathbb{E}\xi_i^2.$
It is easy to see that
for all $x  >0$ and all $0< \varepsilon < 1,$
\begin{eqnarray*}
\mathbb{P} \bigg(\frac{S_n}{[S]_n} \geq x \bigg) &\leq&\mathbb{P} \bigg(\frac{S_n}{[S]_n} \geq x, \     [S]_n \geq B_n^2(1-\varepsilon) \bigg) + \mathbb{P} \bigg(    [S]_n < B_n^2(1-\varepsilon) \bigg)\\
&=&  \mathbb{P} \bigg(\frac{S_n}{[S]_n} \geq x, \     [S]_n \geq B_n^2(1-\varepsilon) \bigg) + \mathbb{P} \bigg(   \sum_{i=1}^n(\xi_i^2- \mathbb{E}\xi_i^2)   < -B_n^2  \varepsilon \bigg).
\end{eqnarray*}
The first term of the last bound can be estimated  by (\ref{ineqPSn02}). For the second term  of the last bound, notice that $(\xi_i^2- \mathbb{E}\xi_i^2)_{i=1,...,n}$ are centered random variables bounded from below, and they are independent once $(\xi_i)_{i=1,...,n}$   are independent. Thus
we need the following  Bernstein type exponential  inequalities for centered random variables bounded from below.

\begin{theorem} \label{thP-x}
Assume that $\xi_i \geq -1$ for all $i \in [1, n]$. Then for all  $x>0$,
\begin{eqnarray}\label{ber01}
\mathbb{P} \bigg(\frac{S_n}{ \langle S \rangle_ n } \leq - x \bigg)
 &\leq&    \inf_{p>1}\bigg(\mathbb{E}\bigg[ \exp\bigg\{-(p-1)  \Big((1+x)\log (1+x) -x \Big)  \langle S \rangle_n  \bigg\} \mathbf{1}_{\{ S_n \leq -x \langle S \rangle_ n   \} }  \bigg]  \bigg)^{  1/ p} \label{ber00} \\
& \leq & \inf_{p>1}\bigg(\mathbb{E}\bigg[ \exp\bigg\{-(p-1) \frac{x^2}{2(1+x/3)}\langle S \rangle_n  \bigg\} \mathbf{1}_{\{ S_n \leq -x \langle S \rangle_ n   \} }  \bigg]  \bigg)^{  1/ p}. \label{ber001}
\end{eqnarray}
\end{theorem}

Inequality (\ref{ber001}) implies that for all  $x>0$,
\begin{eqnarray}
\mathbb{P} \bigg(\frac{S_n}{ \langle S \rangle_ n } \leq - x \bigg)
& \leq & \inf_{p>1}\bigg(\mathbb{E}\bigg[ \exp\bigg\{-(p-1) \frac{x^2}{2(1+x/3)}\langle S \rangle_n  \bigg\}  \bigg]  \bigg)^{ 1/ p}.  \label{ber01}
\end{eqnarray}
It seems that the bound  (\ref{ber01}) is usually decreasing in $p$. For instance, consider the independent case.
When $(\xi_i)_{i=1,\cdots,n}$  are independent random variables, we have $\langle S\rangle _n=\textrm{Var}(S_n),$
where $\textrm{Var}(S_n)$ stands for the variance of $S_n.$
Then the bound (\ref{ber01}) is decreasing in $p$.

For more exponential inequalities similar to that of Theorem  \ref{thP-x}, we refer to Theorem 1.3 of \mbox{de la Pe\~{n}a \cite{VH99}}.  In particular,
de la Pe\~{n}a proved   (\ref{ber001})  with $p=2$. Moreover, de la Pe\~{n}a also proved the following  Bernstein type exponential  inequalities:
  for all $x,  y> 0$,
\begin{eqnarray}
\mathbb{P} \bigg(\frac{S_n}{\langle S\rangle _n} \leq - x, \     \langle S \rangle_n \geq y \bigg) &\leq&   \exp\bigg\{-  \Big((1+x)\log (1+x) -x \Big) y     \bigg\}   \label{ber002} \\
& \leq&   \exp\bigg\{-\frac{x^2 y}{2(1+x/3)}  \bigg\}.  \label{ber02}
\end{eqnarray}
It is easy to see that the inequalities  (\ref{ber002}) and (\ref{ber02}) are respectively the counterparts of (\ref{ineqPSn012}) and (\ref{ineqPSn02}) for $S_n/\langle S\rangle _n$.

Notice that in the independent case, the bounds (\ref{ber01}) and  (\ref{ber02}) are exactly Bernstein's bound (\ref{sdfsdkm}). Thus  (\ref{ber01}) and  (\ref{ber02})  can be  regarded as Bernstein type inequalities for  martingales.

The following deviation inequality for self-normalized martingales has its independent interest.
\begin{theorem}\label{thPbbM1}
Assume that $\xi_i \geq -1$ for all $i \in [1, n]$. Then for all $b> 0,  M\geq 1$ and $x> 0$,
\begin{equation}\label{ineqPbbM1}
\mathbb{P} \bigg(\frac{S_n}{ \sqrt{[S]_ n} } \geq x, \, b \leq   \sqrt{[S]_ n} \leq bM \bigg) \leq   \sqrt{e} \Big(1+2(1+x)\ln M \Big) \exp\bigg\{ -\frac{x^2}{2(1+x/b)}  \bigg\}.
\end{equation}
\end{theorem}

Similarly, when $[S]_n$ in the left hand side of (\ref{ineqPbbM1}) is replaced by $\langle S \rangle_ n$, we have the following inequality for normalized martingales. Such type inequalities are due to Liptser and Spokoiny \cite{Ls01}.
\begin{theorem}\label{thP-xbbM}
Assume that $\xi_i \geq -1$ for all $i \in [1, n]$. Then for all $b> 0,  M\geq 1$ and $x> 0$,
\begin{eqnarray}\label{ineqPbbM2}
\mathbb{P} \bigg(\frac{S_n}{ \sqrt{\langle S \rangle_ n} } \leq -x, \, b \leq   \sqrt{\langle S \rangle_ n} \leq bM \bigg)
 \leq    \sqrt{e} \Big(1+2(1+x)\ln M \Big) \exp\bigg\{ -\frac{x^2}{2(1+x/(3b))}  \bigg\}.
\end{eqnarray}
\end{theorem}

It is interesting to see that in the independent case,   inequality (\ref{ineqPbbM2}) with $b=\sqrt{\textrm{Var}(S_n)}$ and $ M=1$ reduces to exactly Bennett's inequality, up to an absolute constant $\sqrt{e}$. Thus the bound (\ref{ineqPbbM2}) is rather tight.


\section{Applications}

\subsection{Self-scaling sums}
As an application of  Theorem  \ref{thPSn}, we consider the self-normalized sums of i.i.d.\ random variables.
\begin{theorem}\label{the31}
Assume that $(\xi_i)_{i\geq1 }$ are  i.i.d.\ random variables.
Assume that $$\xi_1 \geq -1\ \ \ \   \textrm{and}\ \ \ \    \mathbb{E}\xi_1^{2p} < \infty,$$
  where $1<p\leq 2$.
 Denote $\sigma^2=\mathbb{E}\xi_1^2.$ Then  for all $x >0$ and $y \in (0, \, \sigma^2),$
 \begin{eqnarray*}
 \mathbb{P} \bigg(\frac{S_n}{[S]_n} \geq x  \bigg)   \leq  \exp\bigg\{-\frac{x^2 (\sigma^2 -y)}{2(1+x)}n  \bigg\}+\exp{\bigg\{ -\frac14 \frac{  (p-1)y^{p/(p-1)} }{  (  \mathbb{E}\xi_1^{2p})^{1/(p-1)}}  n \bigg\}}.
\end{eqnarray*}
In particular, it implies that for all $ x   \in (0, \, 1),$
 \begin{eqnarray*}
 \mathbb{P} \bigg(\frac{S_n}{[S]_n} \geq x  \bigg)   \leq  \exp\bigg\{-\frac{ \sigma^2 x^2 }{2(1+ 2x^{(p-1)/p})} n  \bigg\}+\exp{\bigg\{ -  \frac{  (p-1)\, x}{ 4 (  \mathbb{E}\xi_1^{2p}/\sigma^{2p } )^{1/(p-1)}} n\bigg\}}.
\end{eqnarray*}
\end{theorem}

By the last theorem, we have the following moderate deviation result: for any $x>0$ and $\alpha \in ( 0,\,\frac12),$
 \begin{eqnarray*}
\lim_{n \rightarrow \infty} \frac1{n^{1-2\alpha}} \log \mathbb{P} \bigg(\frac{S_n}{[S]_n} \geq  \frac{x}{n^\alpha } \bigg)  \leq  -\frac{1 }{2 } \sigma^2 x^2 .
\end{eqnarray*}
For more such type moderate deviation results, we refer to Shao \cite{S97} and  Jing et al.\ \cite{JLZ12}, 
where the authors established the moderate deviation principles for self-normalized sums $S_n/\sqrt{[S]_n}.$

%

\subsection{Student's t-statistics}
Consider Student's t-statistic $T_n$ defined by
\[
T_n=\sqrt{n} \,  \bar{\xi} \Big/ \Big( \frac{1}{n-1} \sum_{j=1}^n(\xi_j-\bar{\xi})^2 \Big)^{1/2},
\]
where $\bar{\xi}= \sum_{i=1}^n \xi_i/n$. Clearly, $T_n$ and $S_n/ \sqrt{[ S]_n}$ are closely related via the following identity:
\begin{equation}\label{Shao01}
T_n=\frac{S_n}{ \sqrt{[ S]_n} }\bigg( \frac{n-1}{ n-(S_n/ \sqrt{[ S]_n})^2}   \bigg)^{1/2}.
\end{equation}
Since $x/(n-x^2)^{1/2}$ is increasing on $(-\sqrt{n},\sqrt{n})$, it follows from (\ref{Shao01}) that
\begin{equation}\label{Shao02}
\{T_n \geq x \}=\bigg\{ \frac{S_n}{ \sqrt{[ S]_n} } \geq x \Big( \frac{n}{ n+ x^2-1}  \Big)^{1/2} \bigg\}.
\end{equation}
The above fact was pointed out by Efron \cite{E69}.
With the help of (\ref{Shao02}), the following large deviation type result  for t-statistic is an immediate consequence  of Theorem  \ref{thPbbM1}.

\begin{theorem} \label{thPTn}
Assume that $\xi_i \geq -1$ for all $i \in [1, n]$. Then for all $b> 0,  M\geq 1$ and $x> 0$,
\begin{eqnarray}\label{ineqPTn01}
&&\mathbb{P} \big( T_n \geq x, \, b \leq   \sqrt{[S]_ n} \leq bM\big) \nonumber \\
&&\leq   \sqrt{e} \bigg(1+2 x \Big( \frac{n}{ n+ x^2-1}  \Big)^{1/2} \ln M \bigg) \exp\Bigg\{ -\frac{x^2\big( \frac{n}{ n+ x^2-1}  \big)}{2 \big(1 + x \big( \frac{n}{ n+ x^2-1}\big)^{1/2} /b \big)}  \Bigg\}.
\end{eqnarray}
\end{theorem}

\subsection{Autoregressive processes}
The model of autoregressive  process can be expressed as follows: for all $n\geq 0$, by
\begin{equation}\label{eqAuto}
X_{n+1}=\theta X_n+\varepsilon_{n+1}
\end{equation}
where $X_n$ and $\varepsilon_n$ are the observations and driven noises, respectively. We assume that $(\varepsilon_n)$ is a sequence of i.i.d. centered random variables with variation  $\sigma^2 > 0$ and $X_0=\varepsilon_0$.   We can estimate the unknown parameter $\theta$ by the least-squares estimator given by, for all $n\geq 1$,
\begin{equation}\label{eqEstimate}
\hat{\theta}_n=\frac{ \sum_{k=1} ^n X_{k-1} X_k}{ \sum_{k=1} ^n X_{k-1}^2}
\end{equation}
Bercu and Touati  \cite{BT08} has established the convergence rate of $\hat{\theta}_n-\theta$ when $X_0$ and $(\varepsilon_n)$ are   normal random variables.
Here, we would like to give a convergence rate of $\hat{\theta}_n-\theta$ for the case that the driven noises $(\varepsilon_{n})$ are bounded.
Applying Theorem  \ref{thP-x} and de la Pe\~{n}a's inequality (\ref{ber02}),  we have the following exponential inequalities.
\begin{theorem}\label{ththeta}
Assume that  $|\varepsilon_i|  \leq C$ for some positive constant $C$ and all $i.$  If $|\theta| < 1 $, then for all  $x>0$,
\begin{equation}
\mathbb{P} \bigg( |\hat{\theta}_n-\theta|  \geq x  \bigg)
\leq   2\inf_{p>1}\Bigg(\mathbb{E}\Bigg[ \exp\Bigg\{-(p-1)\frac{x^2 }{2\big(\sigma^2 + \frac{ x\, C^2 }{  3(1-|\theta|)} \big)} \sum_{k=1}^n X_{k-1}^2 \Bigg\}  \Bigg]  \Bigg)^{ 1 /p},
\end{equation}
and, for all $ x, y> 0$,
\begin{eqnarray}\label{dfs01}
\mathbb{P} \bigg( |\hat{\theta}_n-\theta| \geq x, \     \sum_{k=1}^n X_{k-1}^2 \geq y \bigg)
 \leq   2  \exp\Bigg\{-\frac{x^2  y}{2\big(\sigma^2+\frac{ x\, C^2 }{  3(1-|\theta|)} \big)}  \Bigg\}.
\end{eqnarray}
\end{theorem}

Inequality (\ref{dfs01}) is similar to an exponential inequalities of de la Pe\~{n}a,   Klass and  Lai \cite{DML07},
which states that when $(\varepsilon_{n})$ are the standard normal random variables, it holds for all $ x, y> 0$,
\begin{eqnarray}
\mathbb{P} \bigg( |\hat{\theta}_n-\theta| \geq x , \     \sum_{k=1}^n X_{k-1}^2 \geq y \bigg)
 \leq   2  \exp\bigg\{-\frac{1}{2 }x^2  y  \bigg\}.
\end{eqnarray}

By Theorem \ref{thP-xbbM}, we obtain the following result.
\begin{theorem}\label{ththeta2}
Assume that  $|\varepsilon_i|  \leq C$ for some positive constant $C$ and all $i.$  If $|\theta| < 1 $, then for all  $b> 0,  M\geq 1$ and $x> 0$,
\begin{eqnarray}
&&\mathbb{P} \bigg(  \big |\widehat{\theta}_n -\theta \big|\sqrt{ \Sigma_{k=1}^n X_{k-1}^2} \geq x , \   b \leq  \sqrt{ \Sigma_{k=1}^n X_{k-1}^2} \leq bM \bigg) \nonumber \\
&&\ \ \ \ \ \ \  \ \ \ \ \ \ \ \leq 2 \sqrt{e} \bigg(1+2(1+\frac{x}{\sigma})\ln M \bigg) \exp\Bigg\{ -\frac{x^2}{2 \big( \sigma^2 +\frac{ x\, C^2 }{  3(1-|\theta|)}\big) }  \Bigg\} .
\end{eqnarray}
\end{theorem}

\section{Proofs of Theorems}
\subsection{Preliminary lemmas}
The following technical lemma is from Fan et al.\ \cite{FGL15}.
For reader's convenience, we shall give a proof following \cite{FGL15}.
\begin{lemma} \label{leESn01}
Assume that $\xi_i \geq -1$ for all $i \in [1, n]$. For any $\lambda \in [0,1)$, denote
\begin{equation*}
U_n(\lambda)=\exp\Big\{ \lambda S_n+ (\lambda+\log(1-\lambda)) [S]_n \Big\}.
\end{equation*}
Then $(  U_i(\lambda),\mathcal{F}_{i})_{i=0,\cdots,n}$  is a   supermartingale,  and satisfies that for all $\lambda \in [0, 1)$,
\begin{equation}\label{ineqleESn}
\mathbb{E}\big[ U_n(\lambda) \big]\leq 1.
\end{equation}
\end{lemma}
Proof.
\vspace{0.3cm}
Assume $\xi_i \geq -1$ and $\lambda \in [0,1)$, then $\lambda \xi_i\geq -\lambda> -1$. We consider the function
\begin{equation*}
f(x)=\frac{\log(1+x)-x}{x^2 /2},\,\,\,\,\,\,\,x>-1,
\end{equation*}
it is increasing in $x$, we obtain that
\begin{eqnarray*}
\log(1+\lambda \xi_i)& \geq & \lambda \xi_i+\frac 1 2(\lambda \xi_i)^2 f(-\lambda)\\
& = & \lambda \xi_i+ \xi_i^2(\lambda+\log(1-\lambda)).
\end{eqnarray*}
Therefore, we have
\begin{equation*}
\exp\big\{ \lambda \xi_i+ \xi_i^2(\lambda+\log(1-\lambda))  \big\}\leq 1+\lambda \xi_i.
\end{equation*}
Since $\mathbb{E} \xi_i = 0$, it follows that
\begin{equation*}
\mathbb{E}\Big[  \exp\big\{ \lambda \xi_i+ (\lambda+\log(1-\lambda)) \xi_i^2  \big\}   \Big] \leq 1.
\end{equation*}
For all $\lambda \in [0,1)$ and $n\geq0$, we have
\[
U_n(\lambda)=U_{n-1}(\lambda)\exp\Big\{\lambda \xi_n+ (\lambda+\log(1-\lambda)) \xi_n^2 \Big\}.
\]
Hence, we deduce that  for all $\lambda \in [0,1)$,
\begin{eqnarray*}
\mathbb{E}[U_n(\lambda)|\mathcal{F}_{n-1}]& = & U_{n-1}(\lambda)\mathbb{E}\Big[ \exp\Big\{\lambda \xi_n+ (\lambda+\log(1-\lambda)) \xi_n^2 \Big\} \Big|\mathcal{F}_{n-1}\Big] \\
&\leq &U_{n-1}(\lambda),
\end{eqnarray*}
which means $(  U_i(\lambda),\mathcal{F}_{i})_{i=0,\cdots,n}$ is a positive supermartingale. Moreover, it holds
\[
\mathbb{E}[U_n(\lambda)]\leq\mathbb{E}[U_{n-1}(\lambda)] \leq ... \leq \mathbb{E}[U_1(\lambda)]\leq1.
\]
%
This completes the proof of Lemma \ref{leESn01}.
\hfill\qed

In order to prove Theorem \ref{thP-x}, we need the following lemma of Freedman \cite{Fr75}.
\begin{lemma}  \label{ineq-Sn}
Assume that $\xi_i \geq -1$ for all $i \in [1, n]$. Denote
\begin{equation*}
W_n(\lambda)=\exp\Big\{- \lambda S_n- (e^\lambda -1-  \lambda ) \langle S\rangle_n \Big\}, \ \ \ \lambda \geq 0.
\end{equation*}
Then  $(  W_i(\lambda),\mathcal{F}_{i})_{i=0,\cdots,n}$ is a   supermartingale,  and satisfies that
\begin{equation}\label{ineqleESn}
\mathbb{E}\big[ W_n(\lambda)\big]  \leq 1.
\end{equation}
\end{lemma}

\subsection{Proof of Theorem \ref{thPSn}}\
We follow the method of Bercu and Touati \cite{BT08}. Let $A_n=\{ S_n \geq x[S]_n\}, x >0.$ By Markov's inequality, H\"{o}lder's inequality and Lemma \ref{leESn01}, we have for all $\lambda \in [0,1)$ and $q>1$,
\begin{eqnarray}\label{ineqPAn01}
\mathbb{P}(A_n)&\leq &\mathbb{E}\bigg[\exp{\bigg\{\frac{\lambda}{q} \Big(S_n -x [S]_n \Big)\bigg\}} \mathbf{1}_{ A_n }  \bigg]\nonumber\\
&=&\mathbb{E}\bigg[ \exp{ \bigg\{ \frac{1}{q} \Big(\lambda S _n +(\lambda+\log(1-\lambda) ) [S]_n \Big) \bigg\}} \exp{ \bigg\{  \frac{1}{q} \Big( -\lambda-\log(1-\lambda)  -\lambda x  \Big) [S]_n \bigg\}  } \mathbf{1}_{ A_n }  \bigg]\nonumber\\
&\leq &\bigg( \mathbb{E}\bigg[\exp{ \bigg\{    \frac{p }{ q } \Big( -\lambda-\log(1-\lambda)  - \lambda x \Big)    [S]_n   \bigg\}} \mathbf{1}_{ A_n } \bigg] \bigg)^{ 1/ p}
\Big( \mathbb{E}\big[ U_n(\lambda)\big] \Big)^{  1/  q }\nonumber\\
&\leq & \bigg( \mathbb{E}\bigg[\exp{ \bigg\{   -\frac{p }{ q } \Big( \lambda+\log(1-\lambda) + \lambda x \Big)    [S]_n  \bigg\}} \mathbf{1}_{ A_n }\bigg] \bigg)^{  1/ p},
\end{eqnarray}
where $p = 1 +p/q $.
Consequently, as $ p/q=p-1$, we can deduce from (\ref{ineqPAn01}) that
\[
\mathbb{P}(A_n)\leq \inf_{p>1}\bigg(\mathbb{E}\bigg[ \exp\Big\{ -(p-1) (\lambda+\log(1-\lambda)+\lambda x ) [S]_n  \Big\} \mathbf{1}_{ A_n } \bigg]  \bigg)^{ 1/ p}.
\]
The right hand side of the last inequality  attains its minimum at $$\underline{\lambda}=\lambda(x)=\frac{x}{1+x},$$ therefore we obtain
\[
\mathbb{P}(A_n)\leq \inf_{p>1}\Bigg(\mathbb{E}\Bigg[ \exp\Big\{-(p-1) ( x-\log(1+x) ) [S]_n  \Big\}  \Bigg] \mathbf{1}_{ A_n } \Bigg)^{  1/ p}.
\]
Using the following inequality
\begin{eqnarray} \label{dsghkl}
x-\log(1+x)
  \geq   \frac{ x^2 }{2( 1+x )}, \ \ \ \ x >0,
\end{eqnarray}
  we  deduce that
\begin{eqnarray*}
\mathbb{P}(A_n)&\leq & \inf_{p>1}\bigg(\mathbb{E}\bigg[ \exp\Big\{-(p-1) ( x-\log(1+x) ) [S]_n  \Big\} \mathbf{1}_{ A_n } \bigg]  \bigg)^{  1/ p}\\
&\leq & \inf_{p>1}\bigg(\mathbb{E}\bigg[ \exp\Big\{-(p-1)\frac{x^2}{2(1+x)} [S]_n \Big\} \mathbf{1}_{ A_n } \bigg]  \bigg)^{ 1 /p},
\end{eqnarray*}
which gives the first two desired inequalities. \hfill\qed

Next we prove the last two desired inequalities.  Denote $B_n=\{S_n \geq x [S]_n ,\, [S]_n  \geq y \}$.
By an argument similar to the proof of (\ref{ineqPAn01}), we deduce that for all  $q>1$,
\begin{eqnarray}
\mathbb{P}(B_n)
&\leq &\bigg( \mathbb{E}\bigg[\exp{ \bigg\{    \frac{p }{ q } \Big( -\underline{\lambda}-\log(1-\underline{\lambda})  - \underline{\lambda} x \Big)    [S]_n   \bigg\}} \mathbf{1}_{ B_n } \bigg] \bigg)^{  1/ p}
\Big( \mathbb{E}\big[ U_n(t)\big] \Big)^{ 1/  q }\nonumber\\
&\leq & \bigg( \mathbb{E}\bigg[\exp{ \bigg\{   -\frac{p }{ q } \Big( \underline{\lambda}+\log(1-\underline{\lambda}) + \underline{\lambda} x \Big)y \bigg\}} \bigg] \bigg)^{  1 /p} \nonumber \\
&=&  \exp\bigg\{- \frac {p-1} p \Big( x-\log(1+x) \Big ) y  \bigg\} .\nonumber
\end{eqnarray}
Therefore, by (\ref{dsghkl}), it holds
\begin{eqnarray}
\mathbb{P}(B_n) &\leq &  \inf_{p>1}\exp\bigg\{- \frac {p-1} p \Big( x-\log(1+x) \Big ) y  \bigg\}  \nonumber  \\
&=&   \exp\bigg\{-   \Big( x-\log(1+x) \Big ) y  \bigg\}    \nonumber \\
& \leq & \exp\bigg\{-  \frac{x^2\, y}{2(1+x)}   \bigg\},  \nonumber
\end{eqnarray}
which gives the last two desired inequalities.

\subsection{Proof of Theorem \ref{thP-x}}
For all $x >0$, denote by $$D_n=\{ -S_n \geq x\langle S \rangle_n\}.$$
By exponential Markov's inequality, we deduce that for all $\lambda \in [0,3)$ and $q>1$,
\begin{eqnarray}
\mathbb{P}(D_n)&\leq &\mathbb{E}\bigg[\exp{\bigg\{\frac{\lambda}{q} \Big(-S_n -x \langle S \rangle_n \Big)\bigg\}} \mathbf{1}_{ D_n }  \bigg]\nonumber\\
&=&\mathbb{E}\bigg[ \exp{ \bigg\{ -\frac{\lambda}{q} S _n -\frac{e^\lambda -1- \lambda  }{ q } \langle S \rangle_n \bigg\}} \exp{ \bigg\{ \Big(\frac{e^\lambda-1-\lambda }{ q }  -\frac{\lambda x}{q} \Big) \langle S \rangle_n \bigg\}  } \mathbf{1}_{ D_n }  \bigg].\nonumber
\end{eqnarray}
Using H\"{o}lder's inequality and Lemma \ref{ineq-Sn}, we have for all $\lambda \in [0,3)$ and $q>1$,
\begin{eqnarray}\label{ineqPDn}
\mathbb{P}(D_n) &\leq &\bigg( \mathbb{E}\bigg[\exp{ \bigg\{ \Big( \frac{p(e^\lambda-1-\lambda)}{ q }  -\frac{p \lambda x}{q} \Big) \langle S \rangle_n   \bigg\}} \mathbf{1}_{ D_n } \bigg] \bigg)^{  1/ p}
\Big( \mathbb{E}\big[ W_n(t)\big] \Big)^{  1/q }\nonumber\\
&\leq & \bigg( \mathbb{E}\bigg[\exp{\bigg\{ \frac p q \Big(e^\lambda-1-\lambda  -\lambda x   \Big) \langle S \rangle_n   \bigg\}} \mathbf{1}_{ D_n }\bigg] \bigg)^{ 1/ p},
\end{eqnarray}
where $p = 1 +p/q $.
Consequently, as $ p/q=p-1$, we can deduce from (\ref{ineqPDn}) that
\begin{eqnarray}\label{ineqdsfPDn}
\mathbb{P}(D_n)\leq \inf_{p>1}\bigg(\mathbb{E}\bigg[ \exp\bigg\{ (p-1) \Big( e^\lambda-1-\lambda-\lambda x \Big) \langle S \rangle_n  \bigg\} \mathbf{1}_{ D_n } \bigg]  \bigg)^{  1/ p}.
\end{eqnarray}
The right hand side of the last inequality attains its minimum at $$\overline{\lambda}=\lambda(x):= \log (1+x)  .$$
Substituting $\overline{\lambda} = \lambda(x)$ in (\ref{ineqdsfPDn}), we obtain
\begin{eqnarray}
\mathbb{P}(D_n) \leq     \inf_{p>1}\bigg(\mathbb{E}\bigg[ \exp\bigg\{-(p-1) \Big((1+x)\log (1+x) -x \Big)    \langle S \rangle_n  \bigg\}\mathbf{1}_{ D_n }  \bigg]  \bigg)^{  1/ p}. \label{thlsf02}
\end{eqnarray}
Using the inequality
\begin{eqnarray} \label{ine2362}
e^{\lambda}-1-\lambda \leq \frac{\lambda^2 }{2(1-\lambda/3  )},\ \ \   \ \lambda \in [ 0, 3) ,
\end{eqnarray} we get  for
all $x\geq0 ,$
\begin{eqnarray}
(1+x)\log (1+x) -x
&=&\inf_{\lambda\geq0}  \Big(  e^\lambda-1-\lambda-\lambda x  \Big)  \nonumber\\
&\leq&\inf_{\lambda\geq0}  \Big(  \frac{\lambda^2 }{2(1-\lambda/3  )} -\lambda x  \Big)  \nonumber\\
&=& \exp\left\{-\frac{x^2}{ 1 +x /3 +   \sqrt{1+2 x/3   } }  \right\} \nonumber\\
&\leq& \exp\left\{-\frac{x^2}{2(1+ x/3 )}\right\},   \nonumber
\end{eqnarray}
where the last inequality follows from the fact $\sqrt{1+2 x/3   } \leq 1+ x/3.  $  Thus, from (\ref{thlsf02}), we obtain for
all $x\geq0 ,$
\begin{eqnarray}
\mathbb{P}(D_n)&\leq&    \inf_{p>1}\bigg(\mathbb{E}\bigg[ \exp\bigg\{-(p-1) \Big((1+x)\log (1+x) -x \Big)    \langle S \rangle_n  \bigg\}  \mathbf{1}_{ D_n }\bigg]  \bigg)^{  1/ p} \nonumber \\
 &\leq&    \inf_{p>1}\bigg(\mathbb{E}\bigg[ \exp\bigg\{-(p-1) \frac{x^2}{ 1+  x/3 +  \sqrt{1+2 x /3}   } \langle S \rangle_n  \bigg\}  \mathbf{1}_{ D_n }\bigg]  \bigg)^{ 1/ p} \nonumber \\
& \leq & \inf_{p>1}\bigg(\mathbb{E}\bigg[ \exp\bigg\{-(p-1) \frac{x^2}{2(1+x/3)}\langle S \rangle_n  \bigg\} \mathbf{1}_{ D_n } \bigg]  \bigg)^{  1 / p}.\nonumber
\end{eqnarray}
This proves (\ref{ber00}) and (\ref{ber01}).
%
\hfill\qed

\subsection{Proof of Theorem \ref{thPbbM1}}

The proof of Theorem \ref{thPbbM1} is based on a modified method of Liptser and Spokoiny \cite{Ls01}. Given $a>1$, introduce the geometric series $b_k=b a^k$ and define random events
$$C_k=\bigg\{ \frac{S_n}{ \sqrt{ [S]_n } } \geq x, \   b_k\leq \sqrt{ [S]_n } <b_{k+1}\bigg\} ,\ \ \  k=0,1,\ldots, K,$$ where $K$ stands for the integer part of $\log_a M$. Clearly, it holds
\begin{eqnarray}
  \bigg\{\frac{S_n}{\sqrt{ [S]_n } } \geq x, \    b\leq \sqrt{ [S]_n } \leq bM \bigg\} \subseteq   \bigcup_{k= 0}^K C_k ,
\end{eqnarray}
which leads to
\begin{eqnarray} \label{ineqCk}
 \mathbb{P} \bigg(\frac{S_n}{\sqrt{ [S]_n } } \geq x, \    b\leq \sqrt{ [S]_n } \leq bM \bigg) \leq \sum_{k=0}^K \mathbb{P}(C_k)  .
\end{eqnarray}
Notice that
$$
 \lambda+\log(1-\lambda)    \geq - \frac{\lambda^2}{2(1-\lambda)}, \ \ \lambda \in [0, 1).
 $$
For any $\lambda \in [0, 1)$, the last inequality and (\ref{ineqleESn}) together implies that

\[
\mathbb{E} \Big[\exp \Big\{ \lambda S_n - \frac{\lambda^2}{2(1-\lambda)} [S]_n  \Big\}   \mathbf{1}_{ C_k }\Big]\leq 1.
\]
Next, taking $\lambda_k=x/(x+b_k),$ for any $x>0$, we obtain
\begin{eqnarray*}
1& \geq &\mathbb{E}\bigg[\exp{ \bigg \{ \frac{x}{x+b_k}S_n - \frac{x^2}{2b_k (x+b_k )} [S]_n \bigg \}} \mathbf{1}_{ C_k }  \bigg]\\
& \geq &\mathbb{E}\bigg[\exp{\bigg\{\frac{x^2}{x+b_k} \sqrt{ [S]_n }- \frac{x^2}{2b_k (x+b_k )} [S]_n \bigg\} } \mathbf{1}_{ C_k }  \bigg]\\
& \geq &\mathbb{E}\bigg[\exp{\bigg\{\inf_{b_k\leq c <  b_{k+1}}\bigg(\frac{x^2 c}{x+b_k}-\frac{x^2 c^2}{2b_k(x+b_k)}\bigg)\bigg\}} \mathbf{1}_{ C_k }  \bigg]\\
& \geq &\mathbb{E}\bigg[\exp{\bigg\{  \frac{x^2 b_{k+1}}{x+b_k}-\frac{x^2 b_{k+1}^2}{2b_k(x+b_k)} \bigg\}} \mathbf{1}_{ C_k }  \bigg],
\end{eqnarray*}
which implies that
\begin{eqnarray*}
\mathbb{P}(C_k)   &\leq&      \exp{\bigg\{  -\frac{x^2  }{ 1+ x/b_k } \Big(a -\frac{a^2} 2 \Big) \bigg\}} \\
 &\leq& \exp{\bigg\{   -\frac{x^2  }{ 1+ x/b  } \Big(a -\frac{a^2} 2 \Big) \bigg\}} .
\end{eqnarray*}
Finally,  we may pick $a$ to make the right-hand side of the last bound possibly small. This leads to the choice $a=1+\frac{1}{1+x},  $ so that
\[
x^2\bigg(a-\frac{a^2}{2}\bigg)= x^2\bigg\{1+ \frac{1}{1+x}-\frac{1}{2}\bigg(1+ \frac{1}{1+x}\bigg)^2\bigg\}\geq\frac{1}{2}(x^2-1).
\]
Since  $\log(1+ \frac{1}{1+x})\geq \frac{1}{2(1+x) }$ for $x\geq 0$, we obtain $\log_a M\leq 2(1+ x)\ln M$ and (\ref{ineqPbbM1}) follows by (\ref{ineqCk}).
\hfill\qed

\subsection{Proof of Theorem \ref{thP-xbbM}}
The proof of  Theorem \ref{thP-xbbM} is similar to
the proof of Theorem \ref{thPbbM1}.
 Given $a>1$, introduce the geometric series $b_k=b a^k$ and define random events
$$H_k=\bigg\{ \frac{-S_n}{ \sqrt{ \langle S \rangle_n } } \geq x, \   b_k\leq \sqrt{ \langle S \rangle_n } <b_{k+1}\bigg\} ,\ \ \  k=0,1,\ldots, K,$$ where $K$ stands for the integer part of $\log_a M$. Clearly, it holds
\begin{eqnarray}
  \bigg\{\frac{-S_n}{\sqrt{ \langle S \rangle_n } } \geq x, \    b\leq \sqrt{ \langle S \rangle_n } \leq bM \bigg\} \subseteq   \bigcup_{k= 0}^K H_k ,
\end{eqnarray}
which leads to
\begin{eqnarray} \label{ineqHk}
 \mathbb{P} \bigg(\frac{-S_n}{\sqrt{ \langle S \rangle_n } } \geq x, \    b\leq \sqrt{ \langle S \rangle_n } \leq bM \bigg) \leq \sum_{k=0}^K \mathbb{P}(H_k)  .
\end{eqnarray}
Notice that
$$
e^\lambda -1-  \lambda   \geq   \frac{\lambda^2}{2(1-\lambda/3)}, \ \ \lambda \in [0, 3).
 $$
For any $\lambda \in [0, 3)$, the last inequality and Lemma  \ref{ineq-Sn} together implies that

\[
\mathbb{E} \Big[\exp \Big\{ \lambda (-S_n) -\frac{\lambda^2}{2(1-\lambda/3)} \langle S \rangle_n  \Big\}   \mathbf{1}_{ H_k }\Big]\leq 1.
\]
Next, taking $\lambda_k=x/(b_k +x /3),$ for any $x>0$, we obtain
\begin{eqnarray*}
1& \geq &\mathbb{E}\bigg[\exp{ \bigg \{ \frac{ x}{b_k + x/3}(-S_n) - \frac{   x^2  }{ b_k (b_k +x/3) } \langle S \rangle_n \bigg \}} \mathbf{1}_{ H_k }  \bigg]\\
& \geq &\mathbb{E}\bigg[\exp{\bigg\{\frac{ x^2 }{b_k+ x /3} \sqrt{ \langle S \rangle_n }- \frac{   x^2  }{ b_k (b_k +x/3) } \langle S \rangle_n \bigg\} } \mathbf{1}_{ H_k }  \bigg]\\
& \geq &\mathbb{E}\bigg[\exp{\bigg\{\inf_{b_k\leq c < b_{k+1}}\bigg(\frac{  x^2 c}{b_k+ x/3}-\frac{   x^2 c^2 }{ b_k (b_k +x/3) }   \bigg)\bigg\}} \mathbf{1}_{ H_k }  \bigg]\\
& \geq &\mathbb{E}\bigg[\exp{\bigg\{  \frac{ x^2 b_{k+1}}{b_k+x /3}-\frac{   x^2  b_{k+1}^2 }{ b_k (b_k +x/3) }   \bigg\}} \mathbf{1}_{ H_k }  \bigg],
\end{eqnarray*}
which implies that
\begin{eqnarray*}
\mathbb{P}(H_k)   &\leq&  \exp{\bigg\{  -\frac{ x^2 }{ 1 + \frac{x}{3b_k} } a + \frac{   x^2   }{   ( 1+\frac{x}{3 b_k } ) }   \frac{a^2} 2 \Big) \bigg\}} \\
&\leq&  \exp{\bigg\{  -\frac{ x^2 }{ 1 + \frac{x}{3b_k} } \Big(a -\frac{a^2} 2 \Big) \bigg\}} \\
&\leq& \exp{\bigg\{   -\frac{x^2  }{ 1 + \frac{x}{3b }  } \Big(a -\frac{a^2} 2 \Big) \bigg\}} .
\end{eqnarray*}
Finally, taking $a=1+\frac{1}{1+x},  $  we obtain  the desired inequality  from (\ref{ineqHk}), with an argument similar to the proof of Theorem \ref{thPbbM1}.
\hfill\qed

\subsection{Proof of Theorem \ref{the31}}
In the proof of Theorem \ref{the31},  we need the following lemma of
Chen et al.\ \cite{Cx16}.
\begin{lemma}\label{lechen}
 Let $(\zeta_i)_{i\geq 1}$  be independent nonnegative random variables with $\mathbb{E} \zeta_i^p < \infty$, where $1<p\leq 2$. Then for all $0< y< \sum_{i=1}^n \mathbb{E} \zeta_i$
\begin{eqnarray}
\mathbb{P} \bigg( \sum_{i=1}^n \zeta_i \leq \sum_{i=1}^n \mathbb{E} \zeta_i -y \bigg) \leq \exp{\bigg\{ -\frac{ (p-1) y^{p/(p-1)} }{4 (\sum_{i=1}^n \mathbb{E} \zeta_i^p)^{1/(p-1)}}  \bigg\}}.
\end{eqnarray}
\end{lemma}

Now, we are in position to prove Theorem \ref{the31}.
For all $x >0$ and $y \in (0, \, \sigma^2),$ we have
\begin{eqnarray}
 \mathbb{P} \bigg(\frac{S_n}{[S]_n} \geq x  \bigg)  &=& \mathbb{P} \bigg(\frac{S_n}{[S]_n} \geq x, \     [S]_n \geq n (\sigma^2 - y) \bigg) + \mathbb{P} \bigg(\frac{S_n}{[S]_n} \geq x, \     [S]_n < n ( \sigma^2 -  y)  \bigg) \nonumber  \\
&\leq & \exp\bigg\{-\frac{x^2 (\sigma^2 -y)}{2(1+x)}n  \bigg\}+\mathbb{P} \bigg(   [S]_n -n\sigma^2 <   - n y   \bigg).  \label{deire01}
\end{eqnarray}
By Lemma \ref{lechen}, we have for all $y \in (0, \, \sigma^2),$
\begin{eqnarray}
  \mathbb{P} \bigg(   [S]_n -n\sigma^2 <   - n y   \bigg) &\leq& \exp{\bigg\{ -\frac{ (p-1) (ny)^{p/(p-1)} }{4 ( n \mathbb{E}\xi_1^{2p})^{1/(p-1)}}  \bigg\}} \nonumber \\
  &=&\exp{\bigg\{ -\frac14 (p-1)   \frac{  y^{p/(p-1)} }{  (  \mathbb{E}\xi_1^{2p})^{1/(p-1)}}  n \bigg\}}.  \label{deire02}
\end{eqnarray}
Combining (\ref{deire01}) and (\ref{deire02}) together, we obtain the first desired inequality. Taking $y=x^{(p-1)/p }\sigma^2 ,$  we obtain the second desired inequality.\hfill\qed

\subsection{Proof of Theorem \ref{ththeta}}
By (\ref{eqAuto}), we have $X_k= \sum_{i=0}^k \theta^{k-i } \varepsilon_i. $ Since $|\theta| < 1$ and $|\xi_i | \leq C, $  we deduce that for all $k,$
$$|X_k| \leq C\sum_{i=0}^k |\theta|^{k-i } \leq \frac{C} {1-|\theta|}.$$
From (\ref{eqAuto}) and (\ref{eqEstimate}), it is easy to see that for all $n \geq 1$,
\begin{equation}\label{eqD1}
\hat{\theta}_n-\theta= \frac{ \sum_{k=1} ^n X_{k-1} \varepsilon_k }{ \sum_{k=1}^n X_{k-1}^2 }.
\end{equation}
For any $i=1,\ldots,n$, set
\[
\xi_i= X_{i-1} \varepsilon_i (1-|\theta|) /C^2\,\,\,\,\textrm{and}\,\,\,\, \mathcal{F}_{i}=\sigma \big(  \varepsilon_k, 0\leq k\leq i \big).
\]
Then $(\xi_i,\mathcal{F}_{i})_{i=1,\ldots,n}$ is a sequence of martingale differences and satisfies
$$|\xi_i| \leq 1  $$
and
\[
 \langle S\rangle_n = \sum_{k=1}^n \mathbb{E} \big[ \xi_i^2  \big|  \mathcal{F}_{i-1} \big]=  \frac{  \sigma^2(1-|\theta|)^2}{ C^4} \sum_{k=1}^n X_{i-1}^2.
\]
Thus we have
\[
\hat{\theta}_n-\theta= \frac{ (1-|\theta|)\sigma^2 }{ C^2  }  \frac{S_n}{ \langle S\rangle _n }.
\]
Applying inequality   (\ref{ber01})  to $(\xi_i,\mathcal{F}_{i})_{i=1,\ldots,n}$,  we deduce that  for all $x  >0,$
\begin{eqnarray}
\mathbb{P} \bigg( \hat{\theta}_n-\theta  \leq -x  \bigg)& =& \mathbb{P} \bigg(\frac{S_n}{\langle S\rangle _n} \leq - \frac{ x \,C^2}{(1- |\theta|)\sigma^2 }  \bigg) \nonumber \\
 &\leq&   \inf_{p>1}\bigg(\mathbb{E}\bigg[ \exp\bigg\{-(p-1)\frac{x^2  }{2(\sigma^2 +x C^2/(3(1-|\theta|)))} \sum_{k=1}^n X_{k-1}^2 \bigg\}\bigg]  \bigg)^{ 1/ p}. \label{thk30}
\end{eqnarray}
Similarly, applying inequality   (\ref{ber01}) to $(-\xi_i,\mathcal{F}_{i})_{i=1,\ldots,n}$, we have  for all $x  >0,$
\begin{eqnarray}
\mathbb{P} \bigg( \hat{\theta}_n-\theta  \geq x  \bigg)
  \leq   \inf_{p>1}\bigg(\mathbb{E}\bigg[ \exp\bigg\{-(p-1)\frac{x^2 }{2( \sigma^2 +x C^2/(3(1-|\theta|)))} \sum_{k=1}^n X_{k-1}^2 \bigg\}\bigg]  \bigg)^{ 1/ p}.  \label{thk31}
\end{eqnarray}
Combining (\ref{thk30}) and (\ref{thk31}) together, we obtain
for all $x  >0,$
\begin{eqnarray*}
 \mathbb{P} \bigg( |\hat{\theta}_n-\theta | \geq x  \bigg)   \leq  2 \inf_{p>1}\bigg(\mathbb{E}\bigg[ \exp\bigg\{-(p-1)\frac{x^2   }{2(\sigma^2 +x C^2/(3(1-|\theta|)))} \sum_{k=1}^n X_{k-1}^2 \bigg\}\bigg]  \bigg)^{ 1/ p},
\end{eqnarray*}
which gives the first desired inequality.
Applying de la Pe\~{n}a's inequality (\ref{ber02}) to $(\xi_i,\mathcal{F}_{i})_{i=1,\ldots,n}$,  we get for all $x, y >0,$
\begin{eqnarray}
\mathbb{P} \bigg( \hat{\theta}_n-\theta  \leq -x , \     \sum_{k=1}^n X_{k-1}^2 \geq y \bigg)& =& \mathbb{P} \bigg(\frac{S_n}{\langle S\rangle _n} \leq - \frac{ x \,C^2}{(1- |\theta|)\sigma^2 }, \     \langle S \rangle_n \geq y \frac{( 1- |\theta|  )^2\sigma^2}{C^4} \bigg) \nonumber \\
 &\leq&    \exp\bigg\{-\frac{x^2  y}{2(\sigma^2 +x C^2/(3(1-|\theta|)))}  \bigg\}. \label{hthk30}
\end{eqnarray}
Similarly, applying de la Pe\~{n}a's inequality (\ref{ber02}) to $(-\xi_i,\mathcal{F}_{i})_{i=1,\ldots,n}$, we have  for all $x, y >0,$
\begin{eqnarray}
\mathbb{P} \bigg(\hat{\theta}_n-\theta \geq x , \     \sum_{k=1}^n X_{k-1}^2 \geq y \bigg)
 \leq     \exp\bigg\{-\frac{x^2 y}{2(\sigma^2 +x C^2/(3(1-|\theta|)))}  \bigg\}.\label{hthk31}
\end{eqnarray}
Combining (\ref{hthk30}) and (\ref{hthk31}) together, we obtain
for all $x, y >0,$
\begin{eqnarray*}
\mathbb{P} \bigg( |\hat{\theta}_n-\theta| \geq x , \     \sum_{k=1}^n X_{k-1}^2 \geq y \bigg)
 \leq   2  \exp\bigg\{-\frac{x^2  y}{2(\sigma^2+x C^2/(3(1-|\theta|)))}  \bigg\},
\end{eqnarray*}
which gives the second desired inequality. \hfill\qed

\subsection{Proof of Theorem \ref{ththeta2}}\
Recall the notations in the proof of Theorem \ref{ththeta}. It is easy to see that
\[
(\widehat{\theta}_n-\theta ) \sqrt{ \Sigma_{k=1}^n X_{k-1}^2 } =  \frac{ \sum_{k=1} ^n X_{k-1} \varepsilon_k }{   \sqrt{ \sum_{k=1}^n X_{k-1}^2 } } = \sigma \frac{S_n}{ \sqrt{ \langle S\rangle _n } }.
\]
Therefore, by Theorem \ref{thP-xbbM}, for all $b> 0,  M\geq 1$ and $x> 0$,
\begin{eqnarray*}
&&\mathbb{P} \bigg(  \big (\widehat{\theta}_n -\theta \big)\sqrt{ \Sigma_{k=1}^n X_{k-1}^2} \leq - x , \   b \leq  \sqrt{ \Sigma_{k=1}^n X_{k-1}^2} \leq bM \bigg) \\
&&\leq  \mathbb{P} \bigg(\frac{S_n}{ \sqrt{\langle S \rangle_ n} } \leq - \frac{x}{\sigma}, \, b  \frac{(1-|\theta|) \sigma}{C^2}  \leq   \sqrt{\langle S \rangle_ n} \leq bM \frac{(1-|\theta|) \sigma}{C^2}\bigg)\\
&&\leq  \sqrt{e} \bigg(1+2(1+\frac{x}{\sigma})\ln M \bigg) \exp\bigg\{ -\frac{x^2}{2( \sigma^2 +xC^2/(3b(1-|\theta|)))}  \bigg\}.
\end{eqnarray*}
Similarly, the same bound holds for the tail probabilities $$\mathbb{P} \bigg(  \big (\widehat{\theta}_n -\theta \big)\sqrt{ \Sigma_{k=1}^n X_{k-1}^2} \geq  x , \   b \leq  \sqrt{ \Sigma_{k=1}^n X_{k-1}^2} \leq bM \bigg).$$
Hence, we have for all $b> 0,  M\geq 1$ and $x> 0$,
\begin{eqnarray*}
&& \mathbb{P} \bigg(  \big |\widehat{\theta}_n -\theta \big|\sqrt{ \Sigma_{k=1}^n X_{k-1}^2}  \geq  x , \   b \leq  \sqrt{ \Sigma_{k=1}^n X_{k-1}^2} \leq bM \bigg) \nonumber\\
 && \ \ \ \ \ \ \ \ \   \ \ \ \ \ \leq 2 \sqrt{e} \bigg(1+2(1+\frac{x}{\sigma})\ln M \bigg) \exp\bigg\{ -\frac{x^2}{2( \sigma^2 +xC^2/(3b(1-|\theta|)))}  \bigg\},
\end{eqnarray*}
which gives the desired inequality.  \hfill\qed

\subsection*{Acknowledgements}
 The work has been partially supported by the
 National Natural Science Foundation of China (Grant nos.\ 11601375 and  11626250).

\end{document}